\documentclass[12pt]{article}
\usepackage{amsthm,amsfonts,amssymb,amscd,amsmath}
\begin{document}

%\begin{flushleft}
%��� 512.7
%\end{flushleft}

\begin{center}
\large \bf Birationally rigid complete intersections \\
of high codimension
\end{center}\vspace{0.5cm}

\centerline{Daniel Evans and Aleksandr Pukhlikov}\vspace{0.5cm}

\parshape=1
3cm 10cm \noindent {\small \quad\quad\quad \quad\quad\quad\quad
\quad\quad\quad {\bf }\newline We prove that a Fano complete
intersection of codimension $k$ and index 1 in the complex
projective space ${\mathbb P}^{M+k}$ for $k\geqslant 20$ and
$M\geqslant 8k\log k$ with at most multi-quadratic singularities
is birationally superrigid. The codimension of the complement to
the set of birationally superrigid complete intersections in the
natural parameter space is shown to be at least $\frac12
(M-5k)(M-6k)$. The proof is based on the techniques of
hypertangent divisors combined with the recently discovered
$4n^2$-inequality for complete intersection singularities.

Bibliography: 23 titles.} \vspace{1cm}

\vspace{1cm}

\noindent 14E05, 14E07\vspace{1cm}

\section*{Introduction}

{\bf 0.1. Complete intersections of index one.} Let $k\geqslant 20$
be a fixed integer. For any integral $k$-uple $\underline{d}=
(d_1,\dots, d_k)$, such that $2\leqslant d_1\leqslant \dots
\leqslant d_k$ set $M=|\underline{d}|-k$, where
$|\underline{d}|=d_1+\cdots +d_k$ and let
$$
{\cal P}(\underline{d})=\prod^k_{i=1}{\cal P}_{d_i,M+k+1}
$$
be the space of $k$-uples of homogeneous polynomials of degree
$d_1,\dots, d_k$, respectively, on the complex projective space
${\mathbb P}={\mathbb P}^{M+k}$. Here the symbol ${\cal P}_{a,N}$
stands for the linear space of homogeneous polynomials of degree
$a$ in $N$ variables which are naturally interpreted as
polynomials on ${\mathbb P}^{N-1}$. We write
$\underline{f}=(f_1,\dots, f_k)\in {\cal P}(\underline{d})$ for an
element of the space ${\cal P}_{\rm }(\underline{d})$. We set also
$$
{\cal P}_{\rm fact}(\underline{d})\subset {\cal P}_{\rm
}(\underline{d})
$$
to be the set of $k$-uples $\underline{f}=(f_1,\dots, f_k)$ such
that the zero set
$$
V(\underline{f})=\{f_1=\cdots =f_k=0\}\subset {\mathbb P}
$$
is an irreducible, reduced and factorial complete intersection of
codimension $k$. Note that for any $\underline{f}\in{\cal P}_{\rm
fact}(\underline{d})$ the projective variety $V(\underline{f})$ is
a primitive Fano variety of index 1, that is,
$$
\mathop{\rm Cl} V(\underline{f})=\mathop{\rm Pic}
V(\underline{f})= {\mathbb Z} H,
$$
where $H$ is the class of a hyperplane section (this is by the
Lefschetz theorem), and $K_{V(\underline{f})}=-H$. Therefore we
may ask if $V=V(\underline{f})$ is birationally rigid or
superrigid (see \cite[Chapter 2]{Pukh13a} for the
definitions).\vspace{0.1cm}

{\bf Theorem 0.1.} {\it Assume that $M\geqslant 8k\log k$. Then
there exists a non-empty Zariski open subset ${\cal P}_{\rm
reg}(\underline{d})\subset {\cal P}_{\rm fact}(\underline{d})$,
such that:\vspace{0.1cm}

{\rm (i)} for every $\underline{f}\in {\cal P}_{\rm
reg}(\underline{d})$ the variety $V=V(\underline{f})$ is
birationally superrigid,\vspace{0.1cm}

{\rm (ii)} the inequality
\begin{equation}\label{05.05.2017.1}
\mathop{\rm codim} (({\cal P}_{\rm }(\underline{d})\setminus {\cal
P}_{\rm reg}(\underline{d}))\subset {\cal P}_{\rm
}(\underline{d}))\geqslant \frac{(M-5k)(M-6k)}{2}
\end{equation}
holds.}\vspace{0.1cm}

Birational superrigidity of generic complete intersections of
index 1 (for $M\geqslant k+7$) was shown in \cite{Pukh01,Pukh14a},
but only non-singular complete intersections were considered
there, so that the complement to the set of birationally
superrigid varieties could be a divisor. In this paper we include
complete intersections with multi-quadratic singularities into
consideration. As a result, we get a much better estimate for the
codimension of the complement: when $k$ is fixed and $M$ is
growing, the codimension is of order $\frac12 M^2$, which is quite
high.\vspace{0.1cm}

We now proceed to explicit definitions of Zariski open subsets in
${\cal P}_{\rm }(\underline{d})$.\vspace{0.3cm}

%%%%%%%%%%%%%%%%%%%%%%%%%%%%%%%%%%%%%%%%%%%%%%%%%%%%%%%%%%%%%%%%%
%%%%%%%%%%%%%%%%%%%%%%%%%%     SUBSECTION 0.2

{\bf 0.2. Complete intersections with multi-quadratic
singularities.} Let us describe the conditions for the
singularities of a complete intersection that guarantee its
factoriality. Take an arbitrary $k$-uple $\underline{f}\in {\cal
P}_{\rm }(\underline{d})$, the zero set $V=V(\underline{f})$ of
which is an irreducible reduced complete intersection of
codimension $k$. Let $o\in V$ be a point. Fix a system of affine
coordinates $(z_1,\dots,z_{M+k})$ on an affine chart ${\mathbb
C}^{M+k}\subset {\mathbb P}$ with the origin at the point $o$.
Write the corresponding dehomogenized polynomials (denoted by the
same symbols) in the form
$$
\begin{array}{l}
f_1=q_{1,1}+q_{1,2}+\cdots +q_{1,d_1},\\ \phantom{f_1}
\dots \\
f_k=q_{k,1}+q_{k,2}+\cdots +q_{k,d_k},
\end{array}
$$
where $q_{i,j}$ is a homogeneous polynomial in $z_*$ of degree
$j$. For a general point $o\in V$
$$
\dim \langle q_{1,1},\dots, q_{k,1}\rangle =k,
$$
that is, $o\in V$ is non-singular. Assume now that $\dim \langle
q_{1,1},\dots, q_{k,1}\rangle\leqslant k-1$, that is to say, $o\in
V$ is a singular point.\vspace{0.1cm}

{\bf Definition 0.1.} The singularity $o\in V$ is a {\it correct
multi-quadratic singularity of type} $2^l$, where $l\in \{1,\dots,
k\}$, if the following conditions are satisfied:\vspace{0.1cm}

\begin{itemize}

\item $\dim \langle q_{1,1},\dots, q_{k,1}\rangle=k-l$,

\item for a general linear subspace $P\subset {\mathbb P}$ of
dimension $\max \{2k+2, k+3l+3\}$, containing the point $o$, the
intersection $V_P=V\cap P$ has an isolated singularity at the
point $o$,

\item for the blow up $\varphi_P\colon V^+_P\to V_P$ of the point
$o$ the exceptional divisor $Q_P=\varphi^{-1}(o)$ is a
non-singular complete intersection of type $2^l$ in the $\max
\{k+l+1,4l+2\}$-dimensional projective space.

\end{itemize}

Note that by Definition 0.1, the codimension of the singular locus
of $V$ near a correct multi-quadratic singularity is at least
$2k+2$.\vspace{0.1cm}

Now let us discuss the conditions of Definition 0.1 in more
detail. There is a subset $I\subset \{1,\dots,k\}$ such that
$|I|=k-l$ and the linear forms $q_{i,1}$, $i\in I$, are linearly
independent:
$$
\langle q_{1,1},\dots, q_{k,1}\rangle=\langle q_{i,1}\,\,|\,\,
i\in I \rangle.
$$
By the genericity of $P$, the restrictions $q_{i,1}|_P$, $i\in I$,
remain linearly independent, so that the zero set
$$
V_{P,I}=\{f_i|_P=0\,\,|\,\, i\in I\}
$$
near the point $o$ is a non-singular complete intersection of
codimension $k-l$. Let
$$
\varphi_{P,I}\colon V^+_{P,I}\to V_{P,I}
$$
be the blow up of the point $o\in V_{P,I}$ with the exceptional
divisor $E_{P,I}=\varphi^{-1}_{P,I}(o)$ being the
$\max\{k+l+1,4l+2\}$-dimensional projective space. Now we can
consider the blow up $\varphi_P$ as the restriction of the blow up
$\varphi_{P,I}$ onto $V_P$, that is, $V^+_P$ is the strict
transform of $V_P$ on $V^+_{P,I}$. In terms of this presentation,
the exceptional divisor $Q_P\subset E_{P,I}$ is given by the set
of $l$ equations
$$
q_{i,2}|_{E_{P,I}}=0, \quad i\in\{1,\dots,k\}\setminus I.
$$
Definition 0.1 requires $Q_P$ to be a non-singular complete
intersection of type $2^l$ in the projective space
$E_{P,I}$.\vspace{0.1cm}

{\bf Definition 0.2.} We say that an irreducible reduced complete
intersection $V=V(\underline{f})$ has at most correct
multi-quadratic singularities if every point $o\in V$ is either
non-singular or a correct multi-quadratic singularity of type
$2^l$ for some $l\in\{1,\dots, k\}$.\vspace{0.1cm}

The set of $k$-uples $\underline{f}\in {\cal P}_{\rm
}(\underline{d})$ such that $V(\underline{f})$ satisfies
Definition 0.2 is denoted ${\cal P}_{\rm mq}(\underline{d})$. The
subset ${\cal P}_{\rm mq}(\underline{d})\subset{\cal P}_{\rm
}(\underline{d})$ is obviously Zariski open. Since for
$\underline{f}\in {\cal P}_{\rm mq}(\underline{d})$ we have
$$
\mathop{\rm codim}\,(\mathop{\rm Sing} V(\underline{f})\subset
V(\underline{f}))\geqslant 2k+2,
$$
by Grothendieck's theorem on parafactoriality of local rings (see
\cite{CL}) the complete intersection $V(\underline{f})$ is a
factorial variety. Therefore, ${\cal P}_{\rm
mq}(\underline{d})\subset {\cal P}_{\rm
fact}(\underline{d})$.\vspace{0.1cm}

{\bf Theorem 0.2.} {\it The following estimate holds:}
\begin{equation}\label{08.05.2017.1}
\mathop{\rm codim} (({\cal P}_{\rm }(\underline{d})\setminus {\cal
P}_{\rm mq}(\underline{d}))\subset {\cal P}_{\rm
}(\underline{d}))\geqslant \frac{(M-4k+1)(M-4k+2)}{2}-(k-1)
\end{equation}

{\bf Remark 0.1.} We construct the subset ${\cal P}_{\rm
reg}(\underline{d})\subset{\cal P}_{\rm mq}(\underline{d})$ below
by removing some additional closed subsets from ${\cal P}_{\rm
mq}(\underline{d})$.\vspace{0.3cm}

%%%%%%%%%%%%%%%%%%%%%%%%%%%%%%%%%%%%%%%%%%%%%%%%%%%%%%%%%%%%%%%%%
%%%%%%%%%%%%%%%%%%%%%%%%%%     SUBSECTION 0.3

{\bf 0.3. Regular complete intersections.} We keep the coordinate
notations of Subsection 0.2 at a point $o\in V$. For brevity and
uniformity we treat the non-singular case $o\not\in\mathop{\rm
Sing} V$ as a multi-quadratic case of type $2^l$ for $l=0$. Let us
place the homogeneous polynomials
$$
q_{i,1},\,\, i\in I, \quad q_{i,j},\,\, j\geqslant 2,
$$
in the {\it standard order}, corresponding to the lexicographic
order of pairs $(i,j)$: $(i_1,j_1)$ precedes $(i_2,j_2)$, if
$j_1<j_2$ or $j_1=j_2$ but $i_1<i_2$. Thus we obtain a sequence
\begin{equation}\label{08.05.2017.2}
h_1, h_2,\dots, h_{M+k-l}
\end{equation}
of $M+k-l$ homogeneous polynomials in $z_*$ of non-decreasing
degrees: $\mathop{\rm deg} h_{e+1}\geqslant \mathop{\rm deg}
h_{e}$.\vspace{0.1cm}

{\bf Definition 0.3.} The point $o\in V$ is {\it regular} if the
sequence of polynomials, which is obtained from
(\ref{08.05.2017.2}) by removing the last $[2\,\log k]-l$ members,
is regular in ${\cal O}_{o,{\mathbb P}}$. (Here $[\cdot]$ means
the integral part of a non-negative real number; if $l>[2\,\log
k]$, we remove no members of the sequence
(\ref{08.05.2017.2}).)\vspace{0.1cm}

In plain words, Definition 0.3 requires that the set of common
zeros of the polynomials $h_e(z)$ in the sequence, obtained from
(\ref{08.05.2017.2}) by removing the last $[2\,\log k]-l$ members,
is of the correct codimension. Since the polynomials $h_*$ are
homogeneous, we may consider them as polynomials on the projective
space ${\mathbb P}^{M+k-1}$ in the homogeneous coordinates
$(z_1:\cdots :z_{M+k})$ and so understand the regularity in the
projective setting.\vspace{0.1cm}

{\bf Definition 0.4.} The complete intersection
$V=V(\underline{f})$, for $\underline{f}\in {\cal P}_{\rm
mq}(\underline{d})$, is {\it regular}, if it is regular at every
point $o\in V$, singular or non-singular. If this is the case, we
write $\underline{f}\in {\cal P}_{\rm
reg}(\underline{d})$.\vspace{0.1cm}

{\bf Theorem 0.3.} {\it Assume that $\underline{f}\in{\cal P}_{\rm
reg}(\underline{d})$. Then $V=V(\underline{f})$ is birationally
superrigid.}\vspace{0.1cm}

{\bf Theorem 0.4.} {\it The following estimate holds:}
\begin{equation}\label{08.05.2017.3}
\mathop{\rm codim} (({\cal P}_{\rm mq}(\underline{d})\setminus
{\cal P}_{\rm reg}(\underline{d}))\subset {\cal P}_{\rm
}(\underline{d}))\geqslant \frac{
(M-5k)(M-6k)}{2}.
\end{equation}

{\bf Proof of Theorem 0.1.} Since the right hand side of
(\ref{08.05.2017.3}) is obviously higher than that of
(\ref{08.05.2017.1}), Theorem 0.1 follows immediately from
Theorems 0.2, 0.3 and 0.4. Q.E.D.\vspace{0.3cm}

%%%%%%%%%%%%%%%%%%%%%%%%%%%%%%%%%%%%%%%%%%%%%%%%%%%%%%%%%%%%%%%%%
%%%%%%%%%%%%%%%%%%%%%%%%%%     SUBSECTION 0.4

{\bf 0.4. The structure of the paper.} Our paper is organized in
the following way. In Section 1 Theorem 0.3 is shown. This is done
by the technique of hypertangent divisors (the constructions can
be found in \cite{Pukh01} or \cite[Chapter 3]{Pukh13a} or
\cite{Pukh2017a}), combined with the recently discovered
$4n^2$-inequality for complete intersection singularities
\cite{Pukh2017b}. We need to take into consideration the fact that
the regularity condition holds, generally speaking, not for the
whole sequence (\ref{08.05.2017.2}), but for a shorter one, so
that the resulting estimates are weaker than in \cite{Pukh01}.
However, we check that they are still sufficient for birational
superrigidity. By the way, the biggest deviation from the
computations in \cite{Pukh01} is for non-singular
points.\vspace{0.1cm}

In Section 2 we show Theorem 0.2. This is rather straightforward
and done by induction on the codimension $k$ of the complete
intersection (here there is no need to assume that $k\geqslant
20$; the case $k=2$ was done in \cite{EvP2017}, $k=1$ in
\cite{EP1}).\vspace{0.1cm}

In Section 3 we show Theorem 0.4. The computations needed for the
proof are really hard; we did our best to make them as clear and
compact as possible. The estimates for the codimension are
obtained by the ``projection'' technique introduced in
\cite{Pukh98b} and also used in \cite{EvP2017}.\vspace{0.3cm}

%%%%%%%%%%%%%%%%%%%%%%%%%%%%%%%%%%%%%%%%%%%%%%%%%%%%%%%%%%%%%%%%%%%
%%%%%%%%%%%%%%%%%%%%%%%   SUBSECTION 0.5

{\bf 0.5. Historical remarks.} The first complete intersection of
codimension at least 2 that was shown to be birationally rigid was
the complete intersection of a quadric and a cubic $V_{2\cdot
3}\subset {\mathbb P}^5$, see \cite{IskPukh96}; for a modernized
exposition, see \cite[Chapter 2]{Pukh13a}. The variety $V_{2\cdot
3}$ was assumed to be general in the sense that it does not
contain lines with ``incorrect'' normal sheaf. Singular complete
intersections $V_{2\cdot 3}\subset {\mathbb P}^5$ were later
studied in \cite{ChG}.\vspace{0.1cm}

Generic complete intersections $V\subset {\mathbb P}^{M+k}$ of
type $\underline{d}$ with $|\underline{d}|=M+k$ and $M\geqslant
2k+1$ were proved to be birationally superrigid in \cite{Pukh01}.
In \cite{Pukh14a} superrigidity was extended to the families with
$M\geqslant k+3$, $M\geqslant 7$ and $d_k=\max\{d_i\}\geqslant 4$,
and in \cite{Pukh13b} to complete intersections of $k_2$ quadrics
and $k_3$ cubics such that $M\geqslant 12$ and $k_3\geqslant 2$.
Today birational superrigidity remains an open problem only for
three infinite series: complete intersections of type
$\underline{d}$, where $\underline{d}$ is
$$
(2,\dots, 2)\quad \mbox{or} \quad (2,\dots, 2,3)\quad
\mbox{or}\quad (2,\dots, 2,4)
$$
and finitely many families with $M\leqslant 11$.\vspace{0.1cm}

In \cite{EP1} a bound for the codimension of the locus of
non-superrigid hypersurfaces of index 1 was given. Such bounds are
important for investigations of birational geometry of Fano fibre
spaces with a higher-dimensional base, see
\cite{Pukh15d,Pukh2017c}. Similar bounds were obtained for
complete intersections of codimension $k=2$ in \cite{EvP2017} and
for double quadrics and cubics (which could be understood as
complete intersections of codimension 2 in a weighted projective
space) in \cite{Johnstone}.\vspace{0.1cm}

For an alternative approach to proving birational superrigidity of
Fano complete intersections in the projective space, see
\cite{Suzuki}. Here are also some other papers on birational
geometry of Fano complete intersections and their generalizations:
[17-22].

%\cite{AhmOkada,ChelWotz2004,Chel2003,Okada1,Okada2,Okada3}.

%%%%%%%%%%%%%%%%%%%%%%%%%%%%%%%%%%%%%%%%%%%%%%%%%%%%%%%%%%%%%%%%%%%%
%%%%%%%%%%%%%%%%%%%%%%%%%%%%%%%%%%%%%%%%%%%%%%%%%%%%%%%%%%%%%%%%%%%%
%%%%%%%%%%%%%%%%%%%%%% SECTION 1

\section{Proof of birational rigidity}

In this section we prove Theorem 0.3. First (Subsection 1.1) we
remind the definition of a maximal singularity and prove that the
centre of a maximal singularity is of codimension at least 3.
After that, in Subsection 1.2 we construct hypertangent divisors.
The construction is standard but singular points need special
attention. In Subsection 1.3 we exclude the case when the centre
of the maximal singularity is not contained in the singular locus
of $V$. In Subsection 1.4 we exclude the case when the centre of
the maximal singularity is contained in the locus of
multi-quadratic points of type $2^l$. Since it follows that a
mobile linear system can not have a maximal singularity, the
variety $V$ is shown to be birationally superrigid.\vspace{0.3cm}

{\bf 1.1. Maximal singularities.} As usual, we prove that a
variety $V=V(\underline{f})$, where $\underline{f}\in {\cal
P}_{\rm reg}(\underline{d})$, is birationally superrigid by
assuming the converse and obtaining a contradiction. So fix a
tuple $\underline{f}\in {\cal P}_{\rm reg}(\underline{d})$ and the
corresponding complete intersection $V=V(\underline{f})$ and
assume that $V$ is not birationally superrigid. This implies
immediately that for some mobile linear system $\Sigma\subset
|nH|$ and some exceptional divisor $E$ over $V$ the Noether-Fano
inequality
$$
\mathop{\rm ord}\nolimits_E\Sigma>n\cdot a(E)
$$
is satisfied, where $a(E)$ is the discrepancy of $E$ with respect
to $V$. In other words, $E$ is a {\it maximal singularity} of
$\Sigma$ (see, for instance \cite[Chapter 2]{Pukh13a}). Let
$B\subset V$ be the centre of $E$ on $V$, an irreducible
subvariety of codimension $\geqslant 2$.\vspace{0.1cm}

{\bf Lemma 1.1.} $\mathop{\rm codim}(B\subset V)\geqslant
3$.\vspace{0.1cm}

{\bf Proof.} Assume the converse: $\mathop{\rm codim}(B\subset
V)=2$. Then $B\not\subset\mathop{\rm Sing}V$, so that the
inequality
$$
\mathop{\rm mult}\nolimits_B\Sigma>n
$$
holds. Consider the self-intersection $Z=(D_1\circ D_2)$ of the
system $\Sigma$, where $D_1,D_2\in \Sigma$ are general divisors.
Obviously, $Z=\beta B+Z_1$, where $\beta>n^2$ and the effective
cycle $Z_1$ of codimension 2 does not contain $B$ as a
component.\vspace{0.1cm}

Let $P\subset{\mathbb P}$ be a general $(2k+1)$-subspace. Since
$\mathop{\rm codim}(\mathop{\rm Sing}V\subset V)\geqslant 2k+2$,
the intersection $V_P=V\cap P$ is non-singular. By Lefschetz, the
numerical Chow group $A^2V_P$ of codimension 2 cycles on $V_P$ is
${\mathbb Z}H^2_P$, where $H_P$ is the class of a hyperplane
section of $V_P$. Setting $Z_P=Z|_P$ and $B_P=B|_P$, we obtain the
inequality
$$
\mathop{\rm deg}\left(Z_P-\beta B_P\right)\geqslant 0.
$$
As $B_P\sim mH^2_P$ for some $m\geqslant 1$, this inequality
implies that
$$
\mathop{\rm deg}V\cdot (n^2-m\beta)\geqslant 0,
$$
which is impossible. Q.E.D. for the lemma.\vspace{0.1cm}

Note that if $\mathop{\rm codim}(B\subset V)\leqslant 2k+1$, then
$B$ is not contained in the singular locus $\mathop{\rm Sing}V$ of
the complete intersection $V$.\vspace{0.3cm}

%%%%%%%%%%%%%%%%%%%%%%%%%%%%%%%%%%%%%%%%%%%%%%%%%%%%%%%%%%%%%%
%%%%%%%%%%%%%%%%%%%%%%%% SUBSECTION 1.2

{\bf 1.2. Hypertangent divisors.} In order to exclude the maximal
singularity $E$, we need the construction of {\it hypertangent
linear systems}. It is well known and published many times (see
\cite{Pukh01} or \cite[Chapter 3]{Pukh13a} or the most recent
application \cite{Pukh2017a}), but some minor modifications are
needed to cover the multi-quadratic case, so we sketch this
construction here. We fix a point $o$ and use the notations of
Subsection 0.2 and work in the affine chart ${\mathbb C}^{M+k}$ of
the space ${\mathbb P}$ with the coordinates $z_1,\dots,z_{M+k}$;
the point $o\in V$ is the origin. Let $j\geqslant 2$ be an
integer. Recall that for some $l\in\{0,1,\dots, k\}$ and a subset
$I\subset \{1,\dots, k\}$, such that $|I|=k-l$, the linear forms
$q_{i,1}$, $i\in I$, are linearly independent, whereas the other
forms $q_{i,1}$, $i\not\in I$, are their linear combinations.
Denote by
$$
f_{i,\alpha}=q_{i,1}+\cdots +q_{i,\alpha}
$$
the truncated $i$-th equation in the tuple
$\underline{f}$.\vspace{0.1cm}

{\bf Definition 1.1.} The linear system
$$
\Lambda(j)= \left\{\left.\left(\sum_{i\in I} q_{i,1}s_{i,j-1}+
\sum^k_{i=1}\sum^{d_i-1}_{\alpha=2}
f_{i,\alpha}s_{i,j-\alpha}\right)\right|_V =0\right\},
$$
where $s_{i,j-\alpha}$ independently run through the set of
homogeneous polynomials of degree $j-\alpha$ in the variables
$z_*$ (if $j-\alpha<0$, then $s_{j-\alpha}=0$), is called the
$j$-{\it th hypertangent system at the point} $o$.\vspace{0.1cm}

For uniformity of notations, we write $\Lambda(1)$ for the {\it
tangent linear system}:
$$
\Lambda(1)= \left\{\left.\left(\sum_{i\in I}
q_{i,1}s_{i,0}\right)\right|_V =0\right\}.
$$
The Zariski tangent space $\{q_{i,1}=0\,|\, i\in I\}$ will be
written as $T$. We set $c(1)=k-l$ and for $j\geqslant 2$
$$
c(j)=k-l+\sharp\{(i,\alpha)\,|\, i=1,\dots,k,\,
1\leqslant\alpha\leqslant\min\{j,d_i-1\}\,\}.
$$
Further, set $m(j)=c(j)-c(j-1)$, where $c(0)=0$, and for
$j=1,\dots,d_k-1$ take $m(j)$ general divisors
$$
D_{j,1},\dots,D_{j,m(j)}
$$
in the linear system $\Lambda(j)$. Putting them into the standard
order, corresponding to the lexicographic order of the pairs
$(j,\alpha)$ (see Subsection 0.3 for a similar procedure), we
obtain a sequence
$$
R_1,\dots,R_{M-l}
$$
of effective divisors on $V$. Set $N_l=M-l$ if $l>[2\,\log k]$ and
$N_l=M-[2\,\log k]$, otherwise. In what follows, we will really
use only the divisors $R_1,\dots, R_{N_l}$, but it is convenient
to keep the entire sequence.\vspace{0.1cm}

{\bf Proposition 1.1.} {\it The equality
$$
\mathop{\rm codim}\nolimits_o \left(\left(\bigcap^{N_l}_{j=1}
|R_j|\right)\subset V\right)=N_l
$$
holds, where $|R_j|$ stands for the support of}
$R_j$.\vspace{0.1cm}

{\bf Proof.} Since
\begin{equation}\label{20.05.2017.1}
f_{i,\alpha}|_V=(-q_{i,\alpha+1}+\dots)|_V
\end{equation}
for $1\leqslant \alpha\leqslant d_i-1$, where the dots stand for
higher order terms in $z_*$, the codimension of the base locus of
the tangent linear system $\Lambda(1)$ near the point $o$ is equal
to $(k-l)$ and of the hypertangent linear system $\Lambda(j)$,
$j\geqslant 2$, is equal to
$$
(k-l)+\mathop{\rm codim}(\{q_{i,\alpha}|_T=0\,|\, 1\leqslant
i\leqslant k,\, 1\leqslant \alpha\leqslant
1+\min\{j,d_i-1\}\,\}\subset T).
$$
Therefore, for a general choice of hypertangent divisors $R_*$,
the equality
$$
\mathop{\rm codim}\nolimits_o \left(\left(\bigcap^{i}_{j=1}
|R_j|\right)\subset V\right)=i
$$
follows from the regularity of the subsequence
$$
h_1,\dots,h_i
$$
of the sequence (\ref{08.05.2017.2}). Now our claim follows
immediately from the regularity condition (see Definition 0.3).
Q.E.D.\vspace{0.1cm}

For a hypertangent divisor $R_i=D_{j,\alpha}$, where $j\in
\{1,\dots, d_k-1\}$ and $\alpha\in\{1,\dots,m(j)\}$, the number
$$
\beta_{l,i}=\beta(R_i)=\frac{j+1}{j}
$$
is its {\it slope}.\vspace{0.1cm}

Set $\varphi\colon V^+\to V$ to be the blow up of the point $o$
with $Q=\varphi^{-1}(o)$ the exceptional divisor. The symbol
$R^+_i$ means the strict transform of $R_i$ on
$V^+$.\vspace{0.1cm}

{\bf Proposition 1.2.} (i) $R^+_i\sim j\varphi^*H-\gamma_i Q$,
{\it where} $\gamma_i\geqslant j+1$.\vspace{0.1cm}

(ii) {\it For any irreducible subvariety $Y\subset V$ of codimension
$\geqslant 2$ such that $Y\not\subset |R_i|$ the algebraic cycle
$(Y\circ R_i)$ of the scheme-theoretic intersection satisfies the
inequality}
$$
\frac{\mathop{\rm mult}\nolimits_o}{\mathop{\rm deg}} (Y\circ
R_i)\geqslant \beta_{l,i}\frac{\mathop{\rm
mult}\nolimits_o}{\mathop{\rm deg}} Y.
$$
(Here the symbol $\mathop{\rm mult}\nolimits_o/\mathop{\rm deg}$
means, as usual, the ratio of the multiplicity at $o$ to the
degree in ${\mathbb P}$.)\vspace{0.1cm}

{\bf Proof.} (i) follows from (\ref{20.05.2017.1}), (ii) follows
from (i). Q.E.D.\vspace{0.3cm}

%%%%%%%%%%%%%%%%%%%%%%%%%%%%%%%%%%%%%%%%%%%%%%%%%%%%%%%%%%%%%%%%
%%%%%%%%%%%%%%%%%%%%%%   SUBSECTION 1.3

{\bf 1.3. The non-singular case.} In the notations of Subsection
1.1, assume that $B\not\subset\mathop{\rm Sing} V$. We want to
show that this case is impossible by obtaining a contradiction. We
write $N$ for $N_0$ and $\beta_i$ for $\beta_{0,i}$ for simplicity
of notations.\vspace{0.1cm}

By \cite[Chapter 2, Section 2]{Pukh13a} the $4n^2$-inequality is
satisfied:
$$
\mathop{\rm mult}\nolimits_B Z>4n^2,
$$
where $Z$ is the self-intersection of the mobile system
$\Sigma\subset |nH|$. Take a point $o\in B$ of general position,
$o\not\in \mathop{\rm Sing} V$, and let $Y_2$ be an irreducible
component of $Z$ with the maximal value of the ratio $\mathop{\rm
mult}\nolimits_o/\mathop{\rm deg}$. Then
$$
\frac{\mathop{\rm mult}\nolimits_o}{\mathop{\rm deg}}
Y_2>\frac{4}{d}.
$$
Take general hypertangent divisors $R_1,\dots, R_M$ as described
in Subsection 1.2. The first $k$ of them, $R_1,\dots, R_k$, are
actually tangent divisors and we know that
$$
\mathop{\rm
codim}\nolimits_o((|R_1|\cap\cdots \cap|R_k|)\subset V)=k.
$$
Proceeding as in Section 1 of \cite{Pukh01}, we construct a
sequence of irreducible subvarieties
$$
Y_2,\dots, Y_k,
$$
such that $\mathop{\rm codim} (Y_i\subset V)=i$, $Y_2$ is an
irreducible component of $Z$ with the maximal value of the ratio
$\mathop{\rm mult}\nolimits_o/\mathop{\rm deg}$, $Y_{i+1}$ is an
irreducible component of the scheme-theoretic intersection
$(Y_i\circ R_{i-1})$ with the maximal value of the ratio
$\mathop{\rm mult}\nolimits_o/\mathop{\rm deg}$ for
$i=2,\dots,k-1$. Therefore, $Y_k\subset V$ is an irreducible
subvariety of codimension $k$, satisfying the inequality
$$
\frac{\mathop{\rm mult}\nolimits_o}{\mathop{\rm deg}}
Y_k>\frac{2^k}{d},
$$
where $d=d_1\cdot\,\cdots\,\cdot d_k=\mathop{\rm deg}
V$.\vspace{0.1cm}

{\bf Lemma 1.2.} $Y_k\not\subset |R_{k-1}|$.\vspace{0.1cm}

{\bf Proof.} Assume the converse: $Y_k\subset |R_{k-1}|$. The
hypertangent divisors being general, this implies that
$$
Y_k\subset \{q_{1,1}|_V= \dots = q_{k,1}|_V=0\}.
$$
However, as $\mathop{\rm codim} (\mathop{\rm Sing}V\subset
V)\geqslant 2k+2$, we can take the section $V_P$ of $V$ by a
generic linear subspace $P\subset {\mathbb P}$ of dimension
$3k+1$, which is a $(2k+1)$-dimensional non-singular complete
intersection in ${\mathbb P}^{3k+1}$. By Lefschetz, the
scheme-theoretic intersection of codimension $k$ on $V_P$
$$
(\{q_{1,1}|_{V_P}=0\}\circ \cdots\circ \{q_{k,1}|_{V_P}=0\})
$$
must be irreducible and reduced. Therefore, the scheme-theoretic
intersection of codimension $k$ on $V$
$$
(\{q_{1,1}|_{V}=0\}\circ \cdots\circ \{q_{k,1}|_{V}=0\})
$$
is irreducible and reduced. By the regularity condition, this
irreducible subvariety has multiplicity precisely $2^k$ at the
point $o$ and the degree $d$. Therefore, it cannot be equal to
$Y_k$. We got a contradiction, proving the lemma.
Q.E.D.\vspace{0.1cm}

By the last lemma, we can proceed in exactly the same way as in
\cite[Section 2]{Pukh01}: form the scheme-theoretic intersection
$(Y_k\circ R_{k-1})$ and obtain an irreducible subvariety
$Y_{k+1}\subset V$ of codimension $k+1$, satisfying the inequality
$$
\frac{\mathop{\rm mult}\nolimits_o}{\mathop{\rm deg}}
Y_{k+1}>\frac{2^{k+1}}{d}.
$$
After that, still following the arguments of \cite[Section
2]{Pukh01}, we use the hypertangent divisors $R_{k+2},\dots, R_N$
to obtain a sequence of irreducible subvarieties $Y_{k+2},\dots,
Y_N$ of codimension $\mathop{\rm codim} (Y_i\subset V)=i$, such
that $Y_i$ is a component of the algebraic cycle $(Y_{i-1}\circ
R_i)$ of the scheme-theoretic intersection of $Y_{i-1}$ and $R_i$
(the regularity condition and genericity of hypertangent divisors
in their linear systems guarantee that $Y_{i-1}\not\subset |R_i|$)
with the maximal value of the ratio $\mathop{\rm
mult}\nolimits_o/\mathop{\rm deg}$. Therefore,
$$
\frac{\mathop{\rm mult}\nolimits_o}{\mathop{\rm deg}} Y_i
\geqslant\beta_i \frac{\mathop{\rm mult}\nolimits_o}{\mathop{\rm
deg}} Y_{i-1}.
$$
The last subvariety $Y_N$ is positive-dimensional and satisfies
the estimate
$$
\frac{\mathop{\rm mult}\nolimits_o}{\mathop{\rm deg}} Y_N>
\gamma=\frac{2^{k+1}}{d}\cdot \prod^N_{i=k+2}\beta_i.
$$

{\bf Proposition 1.3.} {\it The inequality $\gamma\geqslant 1$
holds.}\vspace{0.1cm}

Note that this claim provides the contradiction we need and
excludes the non-singular case.\vspace{0.1cm}

{\bf Proof.} Now it is convenient to use the {\it whole} set
$R_1,\dots, R_M$ of hypertangent divisors, as we have the obvious
identity
$$
d=d_1\cdot\cdots\cdot d_k=\prod^k_{j=1}\prod^{d_j-1}_{\alpha=1}
\frac{\alpha+1}{\alpha}=\prod^M_{i=1}\beta_i.
$$
Recall that $\beta_1=\cdots=\beta_k=2$ and $\beta_{k+1}=\frac32$.
Therefore, $\gamma$ can be re-written as
$$
\gamma=\frac{4}{3d}\prod^N_{i=1}\beta_i=\frac43\beta^{-1},
$$
where
\begin{equation}\label{23.05.2017.1}
\beta=\prod^M_{i=N+1}\beta_i
\end{equation}
and our proposition follows from\vspace{0.1cm}

{\bf Lemma 1.3.} {\it The inequality $\beta<\frac43$
holds.}\vspace{0.1cm}

{\bf Proof of the lemma.} Note first of all that for $j\geqslant
N+1$ we have $\beta_j\leqslant 1+\frac{1}{a}$, where
$a=\left[\frac{M}{k}\right]$. Indeed, assume the converse:
$\beta_{N+1}>1+\frac{1}{a}$. This means that all homogeneous
polynomials $h_{k+1},\dots, h_{k+N}$ in the sequence
(\ref{08.05.2017.2}) are some $q_{i,\alpha}$ with $\alpha<a$.
Therefore,
$$
N\leqslant \sharp \{ q_{i,\alpha}\,|\, 1\leqslant i\leqslant k,\,
2\leqslant\alpha\leqslant a-1\}.
$$
But the right-hand side of this inequality does not exceed
$k\cdot(a-2)<M-k$. So we get:
$$
M-[2\,\log k]<M-k,
$$
which is a contradiction.\vspace{0.1cm}

We have shown that
$$
\beta\leqslant \left(1+\frac{1}{a}\right)^{[2\,\log
k]}\leqslant\left(1+\frac{1}{a}\right)^{a/4}
$$
as $M\geqslant 8k\log k$ by assumption. Therefore,
$\beta<e^{1/4}<\frac43$, as required. Q.E.D. for the
lemma.\vspace{0.1cm}

Proof of Proposition 1.3 is complete.\vspace{0.1cm}

We have shown that the case $B\not\subset\mathop{\rm Sing}V$ is
impossible.\vspace{0.3cm}

%%%%%%%%%%%%%%%%%%%%%%%%%%%%%%%%%%%%%%%%%%%%%%%%%%%%%%%%%%%%%%%%
%%%%%%%%%%%%%%%%%%%%%     SUBSECTION 1.4

{\bf 1.4. The multi-quadratic case.} Assume now that $B$ is
contained in the closure of the locus of multi-quadratic points of
type $2^l$ but not in the closure of the locus of multi-quadratic
points of type $2^j$ for $j\geqslant l+1$. In other words, a
general point $o\in B$ is a singular multi-quadratic point of type
$2^l$. Let us fix this point.\vspace{0.1cm}

{\bf Proposition 1.4.} {\it The self-intersection $Z$ satisfies
the following inequality:} $\mathop{\rm mult}\nolimits_o
Z>2^{l+2}n^2$.\vspace{0.1cm}

{\bf Proof.} This is the $4n^2$-inequality for complete
intersection singularities, see \cite{Pukh2017b}.
Q.E.D.\vspace{0.1cm}

{\bf Remark 1.1.} The condition for a point $o\in V$ to be a
correct multi-quadratic singularity (see Definition 0.1) is in
fact much stronger than what is required in
\cite{Pukh2017b}.\vspace{0.1cm}

Now let us exclude the multi-quadratic case and thus complete the
proof of Theorem 0.3.\vspace{0.1cm}

Assume first that $1\leqslant l\leqslant k-2$. Let
$$
R_1,\dots,R_{k-l}
$$
be general tangent divisors. Since by the regularity condition
$$
\mathop{\rm codim}\nolimits_o \left(\left(\bigcap^{k-l}_{i=1}
|R_i|\right)\subset V\right)=k-l,
$$
we may argue as in the non-singular case and construct a sequence
of irreducible subvarieties
$$
Y_2,\dots,Y_{k-l}
$$
of codimension $\mathop{\rm codim} (Y_i\subset V)=i$, where $Y_2$
is an irreducible component of the cycle $Z$ with the maximal
value of $\mathop{\rm mult}\nolimits_o/\mathop{\rm deg}$ and
$Y_{i+1}$ is an irreducible component of $(Y_i\circ R_{i-1})$ with
the same property. Obviously,
$$
\frac{\mathop{\rm mult}\nolimits_o}{\mathop{\rm deg}}
Y_{k-l}>\frac{2^k}{d}.
$$
By Lefschetz, the scheme-theoretic intersection
$$
(R_1\circ R_2\circ\cdots\circ R_{k-l})
$$
is irreducible and reduced: we make this conclusion, intersecting
that cycle with the section $V_P$ of $V$ by a generic linear
subspace $P$ of dimension $3k+1$, exactly as in the proof of Lemma
1.2 (in fact, in order to apply the Lefschetz theorem, we could
take a subspace $P$ of a smaller dimension here), we conclude that
$Y_{k-l}\not\subset |R_{k-l-1}|$ and construct an irreducible
subvariety $Y_{k-l+1}$, satisfying the inequality
$$
\frac{\mathop{\rm mult}\nolimits_o}{\mathop{\rm deg}} Y_{k-l+1}>
\frac{2^{k+1}}{d}.
$$
After that we argue exactly as in the non-singular case, producing
a sequence of irreducible subvarieties $Y_{k-l+2},\dots, Y_N$, the
last one of which satisfies the estimate
$$
\frac{\mathop{\rm mult}\nolimits_o}{\mathop{\rm deg}}
Y_N>\gamma_l=\frac43 \beta(l)^{-1},
$$
where
\begin{equation}\label{23.05.2017.2}
\beta(l)=\prod^{M-l}_{i=N_l+1}\beta_{l,i}
\end{equation}
(recall that $N_l=M-[2\,\log k]+l$ for $l\leqslant [2\,\log k]$
and $N_l=M-l$, otherwise). The product (\ref{23.05.2017.2})
contains fewer terms than (\ref{23.05.2017.1}) and it is easy to
see that $\beta_{l,M-l-j}=\beta_{M-j}$ for
$j=0,1,\dots,M-l-N_l-1$. Therefore, $\beta(l)<\beta$ for
$l\geqslant 1$ and so $\gamma_l>\gamma>1$, which gives us the
desired contradiction. The multi-quadratic case for $1\leqslant
l\leqslant k-2$ is excluded.\vspace{0.1cm}

Finally, assume that $l\in\{k-1,k\}$. In that case the subvariety
$Y_2$ (an irreducible component of the self-intersection $Z$ with
the maximal value of $\mathop{\rm mult}\nolimits_o/\mathop{\rm
deg}$) satisfies the inequality
$$
\frac{\mathop{\rm mult}\nolimits_o}{\mathop{\rm deg}}
Y_2>\frac{2^{k+1}}{d}
$$
by Proposition 1.4. In this case we omit the part of our arguments
which deals with tangent divisors and proceed straight to the
second part, repeating the arguments for the case $l\leqslant k-2$
word for word.\vspace{0.1cm}

The multi-quadratic case is excluded.\vspace{0.1cm}

Q.E.D. for Theorem 0.3.

%%%%%%%%%%%%%%%%%%%%%%%%%%%%%%%%%%%%%%%%%%%%%%%%%%%%%%%%%%%%%%%%%%%%
%%%%%%%%%%%%%%%%%%%%%%%%%%%%%%%%%%%%%%%%%%%%%%%%%%%%%%%%%%%%%%%%%%%%
%%%%%%%%%%%%%%%%%%%%%% SECTION 2

\section{Irreducible factorial complete intersections}

In this section we prove Theorem 0.2. In Section 2.1 we explain
the strategy of the proof and show the case of a hypersurface.
After that in Section 2.2 we start the inductive part of the
proof, first looking at the easier issue of complete intersections
being irreducible and reduced. Finally, in Subsection 2.3 we
complete the proof, considering complete intersections with
correct multi-quadratic singularities.\vspace{0.3cm}

%%%%%%%%%%%%%%%%%%%%%%%%%%%%%%%%%%%%%%%%%%%%%%%%%%%%%%%%%%%%%%%%
%%%%%%%%%%%%%%%%%%%%%     SUBSECTION 2.1

{\bf 2.1. Complete intersections with correct multi-quadratic
singularities.} Set
$$
\mathcal{P}^{\geqslant j}=\prod_{i=j}^{k}\mathcal{P}_{d_i,M+k+1}
$$
to be the space of truncated tuples $(f_j,\ldots ,f_k)$ and let
$\mathcal{P}_{\rm mq}^{\geqslant j}$ be the set of tuples such
that
$$
V({f_j,\ldots,f_k})=\{f_j=\ldots =f_k=0\}\subset \mathbb{P}
$$
is an irreducible reduced complete intersection of codimension
$k-j+1$ with at most correct multi-quadratic singularities, in the
sense of Definition 0.1 where $k$ is replaced by $k-j+1$. Note
that $\mathcal{P}^{\geqslant 1}=\mathcal{P}(\underline{d})$ and
$\mathcal{P}_{\rm mq}^{\geqslant 1}=\mathcal{P}_{\rm
mq}(\underline{d})$. We will prove Theorem 0.2 by decreasing
induction on $j=k,k-1,\ldots ,1$ in the following form
\begin{equation}\label{est1}
\mathop{\rm codim}((\mathcal{P}^{\geqslant j}\setminus
\mathcal{P}_{\rm mq}^{\geqslant j}) \subset
\mathcal{P}^{\geqslant j})\geqslant
\frac{(M-4k+1)(M-4k +2)}{2}-(k-1).
\end{equation}
The basis of the induction is the case of a hypersurface
$V(f_k)\subset\mathbb{P}$ of degree $d_k$. It is easy to calculate
that the closed subset of reducible or non-reduced polynomials of
degree $d_k$ has codimension
$$
{{M+k+d_k-1}\choose{d_k}}-(M+k+1)
$$
in $\mathcal{P}_{d_k,M+k+1}$ (which corresponds to the case when
$f_k$ has a linear factor), and the closed subset of polynomials
$f_k$ such that the hypersurface $V(f_k)$ has at least one
singular point, which is not a quadratic singularity of rank at
least $7$ has codimenison
$$
\frac{(M+k-6)(M+k-5)}{2}+1
$$
in $\mathcal{P}_{d_k,M+k+1}$ (see a similar detailed calculation
in \cite{EP1} for the case of rank at least 5). Therefore, the
inequality (\ref{est1}) is true for $k=1$. \vspace{0.1cm}

Now let us proceed to the inductive argument.\vspace{0.3cm}

%%%%%%%%%%%%%%%%%%%%%%%%%%%%%%%%%%%%%%%%%%%%%%%%%%%%%%%%%%%%%%%%
%%%%%%%%%%%%%%%%%%%%%%%%   subsection 2.2

{\bf 2.2. The step of induction: irreducibility.} We assume that
(\ref{est1}) is shown for $j+1$. The task is, for a fixed tuple
$(f_{j+1},\ldots ,f_k)\in \mathcal{P}_{\rm mq}^{\geqslant j+1}$,
to estimate the codimension of the set of polynomials $f_j\in
\mathcal{P}_{d_j,M+k+1}$ such that $V(f_j,\ldots ,f_k)$ does not
satisfy the required condition, that is, $(f_j, f_{j+1},\ldots
f_k)\notin \mathcal{P}_{\rm mq}^{\geqslant j}$.\vspace{0.1cm}

Let us first consider the issue of irreducibility and reducedness.
Since by the inductive assumption and the Grothendieck theorem
\cite{CL} $V(f_{j+1},\ldots,f_k)$ is a factorial complete
intersection, we have the isomorphism
$$
{\rm Cl\hspace{0.1cm}} V(f_{j+1}, \ldots ,f_k)\cong
{\rm Pic \hspace{0.1cm}}  V(f_{j+1}, \ldots ,f_k)\cong \mathbb{Z}H,
$$
where $H$ is the class of a hyperplane section, and moreover, for
every $a\in \mathbb{Z}_+$ the restriction map
$$
r_a :{\rm H}^0 (\mathbb{P}, \mathcal{O}_{\mathbb{P}}(a))\rightarrow
{\rm H}^0 (V_{j+1}, \mathcal{O}_{V_{j+1}}(a))
$$
is surjective (where for simplicity of notation we write $V_{j+1}$
for $V(f_{j+1},\ldots ,f_k)$). For $a < d_{j+1}$ it is also
injective, and for $a=d_{j+1}$ we have
$$
{\rm dim\hspace{0.1cm} Ker}\hspace{0.1cm} r_a
= \# \{i\in\{j+1,\ldots ,k\}\hspace{0.1cm} |\hspace{0.1cm}
d_i=d_{j+1}\}.
$$
Now easy calculations show that the set of polynomials $f_j \in
\mathcal{P}_{d_j,M+k+1}$ such that $V(f_j,f_{j+1},\ldots ,f_k)$ is
either reducible or non-reduced, is of codimension at least
$$
{{M+k+d_j-1}\choose{d_j}}-(M+k+1)-(k-j)
$$
(again, this corresponds to the case when the divisor
$\{f_j|_{V_{j+1}}=0\}$ has a hyperplane section of $V_{j+1}$ as a
component). This estimate is higher (and, in fact, much higher)
than what we need so we may assume that $V(f_j,f_{j+1},\ldots
,f_k)$ is irreducible and reduced.\vspace{0.1cm}

Finally, we need to consider the condition for the singularities
of the complete intersection $V(f_j,f_{j+1},\ldots ,f_k)$ to be
multi-quadratic. In order to avoid cumbersome formulae, we will
consider the final case $j=1$ only, when the estimate is the
weakest. For higher values of $j$ the arguments are identically
the same, just the indices and dimensions need to be adjusted
appropriately.\vspace{0.3cm}

%%%%%%%%%%%%%%%%%%%%%%%%%%%%%%%%%%%%%%%%%%%%%%%%%%%%%%%%%%%%%%
%%%%%%%%%%%%%%%%%%%%%%%%%%%%%%%%%%%%%%%%%%%%%%%%%%%%%%%%%%%%%%

{\bf 2.3. Multi-quadratic singularities.} Fix a point $o\in
{\mathbb P}$ and consider a tuple
$({f_1,\ldots,f_k})\in\mathcal{P}^{\geqslant 1}$ with $o\in
V=V({f_1,\ldots,f_k})$. Fix a system of affine coordinates
$(z_1,\ldots ,z_{M+k})$ on an affine chart
$\mathbb{C}^{M+k}\subset \mathbb{P}$ with the origin at the point
$o$. Write the corresponding dehomogenized polynomials (denoted by
the same symbols) in the form
$$
\begin{array}{l}
f_1=q_{1,1}+q_{1,2}+\cdots +q_{1,d_1},\\ \phantom{f_1}
\dots \\
f_k=q_{k,1}+q_{k,2}+\cdots +q_{k,d_k},
\end{array}
$$
where $q_{i,j}$ is a homogeneous polynomial in $z_*$ of degree
$j$. Assume that
$$
\dim\langle q_{1,1},\dots, q_{k,1}\rangle =k-l,
$$
with $l\geqslant 0$. Let $I \subset \{1,\ldots ,k\}$ be a subset
with $|I|=k-l$ such that the linear forms $\{q_{i,1}\hspace{0.1cm}
| \hspace{0.1cm} i\in I\}$ are linearly independent. Set $\Pi
\subset \mathbb{C}^{M+k}$ to be the subspace
$$
\Pi =\{q_{i,1}=0\hspace{0.1cm} | \hspace{0.1cm} i\in I\}\cong {\mathbb C}^{M+l}.
$$
By assumption, for every $j\in J=\{1,\dots, k\}\setminus I$ there
are (uniquely determined) constants $\beta_{j,i}$, $i\in I$, such
that
$$
q_{j,1}=\sum_{i\in I}\beta_{j,i}q_{i,1}.
$$
Set for every $j\in J$
$$
q^*_{j,2}=\left.\left( q_{j,2}-
\sum_{i\in I}\beta_{j,i}q_{i,2}\right)\right|_{\Pi}.
$$
The following statement translates the condition for the point $o$
to be a correct multi-quadratic singularity into the language of
properties of the quadratic forms $q^*_{j,2}$ introduced
above.\vspace{0.1cm}

{\bf Proposition 2.1.} {\it Assume that for a general subspace
$\Theta\subset{\mathbb P}(\Pi)$ of dimension
$$
b=\mathop{\rm max} \{k+l+1, 4l+2\}
$$
the set of quadratic equations
$$
\left\{ q^*_{j,2}|_{\Theta}=0\,\, |\,\, j\in J\right\}
$$
defines a non-singular complete intersection of type $2^l$. Then
$o\in V$ is a correct multi-quadratic singularity of type
$2^l$.}\vspace{0.1cm}

{\bf Proof.} Indeed, it is easy to see that the germ $o\in V$ is
analytically equivalent to the closed set in $\Pi$ defined by $l$
equations
$$
0=q^*_{j,2}+\dots,\quad j\in J,
$$
where the dots stand for higher order terms. The rest is obvious.
Q.E.D.\vspace{0.1cm}

{\bf Remark 2.1.} In the notations of Definition 0.1, the
exceptional divisor $Q_P$ is precisely the complete intersection
of $l$ quadrics $\{q^*_{j,2}|_{\Theta}=0\}$, $j\in J$, in the
$b$-dimensional space $\Theta$. Proposition 2.1 gives a sufficient
condition for the point $o$ to be a correct multi-quadratic
singularity. Now we use this criterion to estimate the codimension
of the set of tuples violating the conditions of Definition 0.1 at
the given point $o\in V$.\vspace{0.1cm}

{\bf Definition 2.1.} We say that an $l$-uple $(q^*_{j,2}\,\,
|\,\, j\in J)$ is {\it correct}, if its zero set in ${\mathbb
P}(\Pi)$ is an irreducible reduced complete intersection $Q_{\Pi}$
satisfying the inequality
$$
\mathop{\rm codim}(\mathop{\rm Sing} Q_{\Pi}\subset
Q_{\Pi})\geqslant b.
$$

{\bf Corollary 2.1.} {\it Assume that the $l$-uple $(q^*_{j,2}\,\,
|\,\, j\in J)$ is correct. Then $o\in V$ is a correct
multi-quadratic singularity of type $2^l$.}\vspace{0.1cm}

Since in the subsequent arguments (up to the end of this section)
only the quadratic forms $q_{i,2}$ will be involved, we may assume
without loss of generality that
$$
J=\{1,\dots, l\}
$$
and $I=\{l+1,\dots, k\}$. Fixing the forms $q_{i,2}$ for $i\in I$,
we work with the $l$-uples
$$
(q^*_{j,2}\,|\, j=1,\dots, l)\in {\cal P}^{\times l}_{2, M+l}.
$$
Theorem 0.2 is obviously implied by the following
proposition.\vspace{0.1cm}

{\bf Proposition 2.2.} {\it The codimension of the closed set
${\cal X}\subset {\cal P}^{\times l}_{2, M+l}$ of incorrect
$l$-uples is at least}
\begin{equation}\label{27.02.2018.1}
\frac{(M+3-b)(M+4-b)}{2}-(l-1).
\end{equation}

(Recall that $b=\mathop{\rm max} \{k+l+1, 4l+2\}$.)\vspace{0.1cm}

{\bf Proof.} Elementary computations show that the codimension of
the closed subset ${\cal X}_*\subset {\cal P}^{\times l}_{2, M+l}$
of linearly dependent $l$-uples is higher than
(\ref{27.02.2018.1}), so we may assume the forms $q^*_{j,2}$,
$j=1,\dots, l$, to be linearly independent. The symbol $Q_{\Pi}$
stands for their zero set. By the symbol $\mathop{\rm Sing}
Q_{\Pi}$ we denote the closed set of points $p\in Q_{\Pi}$, such
that the linear terms of dehomogenised polynomials $q^*_{j,2}$
with respect to any system of affine coordinates with the origin
at $p$ are not linearly independent. (We argue in this way in
order to avoid a discussion of the zero scheme of the forms
$q^*_{j,2}$, $j=1,\dots, l$, being irreducible and reduced at this
stage of the proof.) For
$$
\lambda=(\lambda_1 : \cdots :\lambda_l)\in {\mathbb P}^{l-1}
$$
set
$$
W(\lambda)=\{\lambda_1q^*_{1,2}+\cdots
+\lambda_lq^*_{l,2}=0\}\subset {\mathbb P}^{M+l-1}
$$
to be the corresponding quadric hypersurface in the linear system
generated by $(q^*_{j,2})$. We will use the following simple
observation, which for $k=2$ was used in
\cite{EvP2017}.\vspace{0.1cm}

{\bf Lemma 2.1.} {\it For any point $p\in\mathop{\rm Sing}
Q_{\Pi}$ there is $\lambda\in {\mathbb P}^{l-1}$ such that
$p\in\mathop{\rm Sing} W(\lambda)$.}\vspace{0.1cm}

{\bf Proof.} Obvious computations. Q.E.D. for the
lemma.\vspace{0.1cm}

{\bf Corollary 2.2.} {\it The following inclusion holds:}
$$
\mathop{\rm Sing}
Q_{\Pi}\subset\mathop{\bigcup}\limits_{\lambda\in {\mathbb
P}^{l-1}}\mathop{\rm Sing} W(\lambda).
$$

Set ${\cal R}_{\leqslant a} \subset {\cal P}_{2, M+l}$ to be the
closed subset of quadratic forms of rank $\leqslant a$. It is well
known that
$$
\mathop{\rm codim} ({\cal R}_{\leqslant a} \subset {\cal P}_{2,
M+l})=\frac{(M+l+1-a)(M+l+2-a)}{2}.
$$
Now for every $e=1,\dots, l$ consider the closed subset ${\cal
X}_{e,a} \subset {\cal P}^{\times e}_{2, M+l}$, consisting of
$e$-uples $(g_1,\dots, g_e)$ such that the linear span $\langle
g_1,\dots, g_e\rangle$ has a positive-dimensional intersection
with ${\cal R}_{\leqslant a}$.\vspace{0.1cm}

{\bf Lemma 2.2.} {\it The following estimate holds:}
$$
\mathop{\rm codim} ({\cal X}_{e,a} \subset {\cal P}^{\times e}_{2,
M+l})\geqslant \mathop{\rm codim}({\cal R}_{\leqslant a} \subset
{\cal P}_{2, M+l}) - (e-1).
$$

{\bf Proof.} Consider the natural projections of ${\cal P}^{\times
e}_{2,M+l}= {\cal P}^{\times (e-1)}_{2,M+l}\times {\cal
P}^{}_{2,M+l}$ onto the last factor and the direct product ${\cal
P}^{\times (e-1)}_{2,M+l}$ of the first $e-1$
factors.\vspace{0.1cm}

For any tuple
$$
(g_1,\dots,g_e)\in {\cal P}^{\times e}_{2,M+l}
$$
such that
$$
(g_1,\dots,g_{e-1})\not\in {\cal X}_{e-1,a}
$$
the condition $(g_1,\dots,g_{e})\in {\cal X}_{e,a}$ implies that
the quadratic form $q_e$ belongs to the cone with the base ${\cal
R}_{\leqslant a}$ and the vertex space $\langle g_1,\dots,
g_{e-1}\rangle$, which has dimension at most $\mathop{\rm dim}
{\cal R}_{\leqslant a}+(e-1)$. Arguing by increasing induction on
$e=1,\dots, l$, we complete the proof. Q.E.D. for the
lemma.\vspace{0.1cm}

Now we can complete the proof of Proposition 2.2. Let us consider
an $l$-uple $(q^*_{j,2}\,|\, j=1,\dots, l)$ such that
$$
\mathop{\rm codim} (\mathop{\rm Sing} Q_{\Pi}\subset
Q_{\Pi})\leqslant b-1
$$
or, equivalently, that
$$
\mathop{\rm dim} (\mathop{\rm Sing} Q_{\Pi})\geqslant M+l-b.
$$
By Corollary 2.2 we conclude that the inequality
$$
\mathop{\rm max}\limits_{\lambda\in {\mathbb
P}^{l-1}}\{\mathop{\rm dim} \mathop{\rm Sing}
W(\lambda)\}\geqslant M+1-b
$$
is satisfied, which, in its turn, implies that
$$
(g_1,\dots,g_{e})\not\in {\cal X}_{l,a}
$$
for $a=l+b-2$. Now by Lemma 2.2 we conclude that in the proof of
Proposition 2.2 we can consider only $l$-uples satisfying the
inequality
\begin{equation}\label{27.02.2018.2}
\mathop{\rm codim} (\mathop{\rm Sing} Q_{\Pi}\subset
Q_{\Pi})\geqslant b.
\end{equation}
The rest is very easy. If $Q_{\Pi}$ is an irreducible reduced
complete intersection, then (\ref{27.02.2018.2}) guarantees that
the tuple of quadratic forms under consideration is correct.
Moreover, if for some $e\geqslant 1$ the system of quadratic
equations
$$
q^*_{1,2}=\cdots =q^*_{e,2}=0
$$
defines an irreducible reduced complete intersection of $e$
quadrics, then by (\ref{27.02.2018.2}) it is factorial. Now
arguing as in Subsection 2.2, we can estimate the codimension of
the set of tuples, the zero set of which is not an irreducible
reduced complete intersections. It is easy to check that the
codimension is equal to
$$
\frac{(M+l-1)(M+l-2)}{2}-e.
$$
This completes the proof of Proposition 2.2. Q.E.D.

This completes the proof of Theorem 0.2 as well, as the minimum of
the estimate obtained in Proposition 2.2 occurs for $l=k$.

%%%%%%%%%%%%%%%%%%%%%%%%%%%%%%%%%%%%%%%%%%%%%%%%%%%%%%%%%%%%%%%%%%%%
%%%%%%%%%%%%%%%%%%%%%%%%%%%%%%%%%%%%%%%%%%%%%%%%%%%%%%%%%%%%%%%%%%%%
%%%%%%%%%%%%%%%%%%%%%% SECTION 3

\section{Regular complete intersections}

In this section we prove Theorem 0.4. In Subsection 3.1 we produce
the estimates for the codimension of the set of non-regular tuples
of polynomials, given by the projection method. After that, the
proof of Theorem 0.4 is reduced to showing a purely analytical
fact: estimating the minimum of an integral sequence, consisting
of certain binomial coefficients, depending on several integral
parameters. The required computations are quite non-trivial. We
perform them in several steps. In Subsection 3.2 a number of
reductions simplifies the task. In Subsection 3.3 we employ the
classical Stirling formula to approximate with good precision the
expressions to be minimized by a smooth function and study that
function using the standard tools of calculus. In Subsections 3.4
and 3.5 we complete the proof, showing the required
estimates.\vspace{0.3cm}

{\bf 3.1. The projection method.} We use the notations of
Subsection 0.3. Since an elementary dimension count relates the
codimension of the set of globally non-regular tuples
$\underline{f}$ (which is what Theorem 0.4 estimates) to the
codimension of the set of tuples $\underline{f}$ that are
non-regular at a fixed point $o\in V(\underline{f})$ (see Theorem
3.1 and the comments below), we concentrate on the local problem:
fix a point $o\in {\mathbb P}$, a system of affine coordinates
$z_1,\dots, z_{M+k}$ with the origin at $o$ and consider
(non-homogeneous) tuples $\underline{f}$ such that $o\in
V(\underline{f})$.\vspace{0.1cm}

Next, we fix $l\in \{0,1,\dots, k\}$ and assume that the rank of
the set of linear forms $q_{i,1}$, $i=1,\dots, k$, is equal to
$k-l$, so that in the sequence (\ref{08.05.2017.2}) exactly the
first $k-l$ polynomials are linear forms. We fix them, either, so
that the linear subspace
$$
\Pi = \{h_1=\ldots =h_{k-l}=0\} \cong \mathbb{C}^{M+l}
$$
of the space ${\mathbb C}^{M+k}_{z_*}$ is also fixed. Recall the
notation
$$
N_l=M-\mathop{\rm max} \{[2\,\mathop{\rm log} k],l\},
$$
introduced in Subsection 1.2. Set
$$
g_i=h_{k-l+i}|_{{\mathbb P}(\Pi)},
$$
$i=1,\dots, N_l$. This is a sequence of $N_l$ homogeneous
polynomials of non-decreasing degrees $m_i=\mathop{\rm deg} g_i$
on the projective space ${\mathbb P}(\Pi)\cong {\mathbb
P}^{M+l-1}$. Define the space of such sequences:
$$
{\cal G}(\underline{d},l)=\prod^{N_l}_{i=1} {\cal P}_{m_i, M+l}.
$$
It is obvious that the point $o\in V$ is regular (as a
multi-quadratic point of type $2^l$ in the sense of Definition
0.3) if and only if the sequence
$$
g_1,\dots, g_{N_l}
$$
is regular, that is to say, if the closed algebraic set
$$
\{g_1=\cdots =g_{N_l}=0\}\subset {\mathbb P}(\Pi)
$$
has codimension $N_l$. Set ${\cal Y}={\cal
Y}(\underline{d},l)\subset {\cal G}(\underline{d},l)$ to be the
closed set of {\it non}-regular tuples.\vspace{0.1cm}

{\bf Theorem 3.1.} {\it Assume that $M \geqslant 8k \log k$ and $k
\geqslant 20$. Then}
$$
\mathop{\rm codim} ({\cal Y}\subset {\cal
G}(\underline{d},l))\geqslant \frac{(M-5k)(M-6k)}{2}+M+k.
$$

Taking into account that the point $o$ varies in the
$M+k$-dimensional projective space ${\mathbb P}$ and the original
tuple $\underline{f}$ satisfies the conditions $f_1(o)=\cdots
=f_k(o)=0$ and $\dim \langle q_{i,1}\,|\, 1\leqslant i\leqslant
k\rangle=k-l$, an elementary dimension count gives Theorem 0.4 as
an immediate corollary of Theorem 3.1.\vspace{0.1cm}

The rest of this section is the {\bf proof of Theorem 3.1}. Our
main tool is {\it the projection method}, developed in
\cite{Pukh98b} and explained and applied to solving similar
problems in \cite[Chapter 3]{Pukh13a} and many papers, e.g.
\cite{EvP2017,EP1}. The idea is to represent
$$
{\cal Y}=\coprod^{N_l}_{e=1} {\cal Y}_e
$$
as a disjoint union of constructive subsets ${\cal Y}_e$,
consisting of tuples $(g_1,\dots, g_{N_l})$ such that the closed
set
$$
\{g_1=\cdots =g_{e-1}=0\}\subset {\mathbb P}(\Pi)
$$
is of codimension $e-1$, but $g_e$ vanishes on some irreducible
component of that set (if $e=1$, this means simply that the
quadratic form $g_1$ is identically zero). The projection method
estimates the codimension of ${\cal Y}_e$ in ${\cal
G}(\underline{d},l)$ as follows:
$$
\mathop{\rm codim}({\cal Y}_e\subset
{\cal G}(\underline{d},l)\geqslant \gamma(e,\underline{d},l)=
h^0(\mathbb{P}^{M+l-e},\mathcal{O}_{\mathbb{P}^{M+l-e}}(m_e))=
{M+l-e+m_e \choose M+l-e},
$$
where $m_e=\mathop{\rm deg} g_e$, see, for instance, \cite[Chapter
3]{Pukh13a}. Therefore, in order to prove Theorem 3.1, we must
show that the numbers $\gamma(e,\underline{d},l)$ for $e=1,\dots,
N_l$ are not smaller than the right hand side of the inequality of
Theorem 3.1. This is what we are going to do. The task is quite
non-trivial. First, we do some preparatory work in order to
simplify the inequalities to be shown and reduce the number of
integral parameters, on which the numbers
$\gamma(e,\underline{d},l)$ depend.\vspace{0.3cm}

%%%%%%%%%%%%%%%%%%%%%%%%%%%%%%%%%%%%%%%%%%%%%%%%%%%%%%%%%%%%%%%%%%%%%%%%%%%%
%%%%%%%%%%%%%%%%%%%%%%%%   subsection 3.2.

{\bf 3.2. Reductions.} If the original tuple $\underline{f}$ of
defining polynomials consists of $k_2$ quadrics, $k_3$ cubics,
$\ldots$, $k_m$ polynomials of degree $m=d_k \geqslant 8\log k$,
then
$$
k_2+k_3+\dots + k_m =k
$$
and
$$
2k_2+3k_3+\dots +mk_m = |\underline{d}|=d_1+\dots +d_k.
$$
It is easy to see that
$$
m_e=\mathop{\rm deg}g_e=
\mathop{\rm min}\left\{ j\, \left| \, \sum_{\alpha=2}^{j}
\left( \sum_{\beta =\alpha}^{m}k_{\beta}\right)\geqslant
e\right\}\right. .
$$
This explicit presentation gives us the first
reduction.\vspace{0.1cm}

{\bf Proposition 3.1.} {\it The following estimate holds:
$$
\gamma(e,\underline{d},l)\geqslant \gamma(e,\underline{d}^*,l),
$$
where the $k$-uple $\underline{d}^*=(d_1^*,\ldots , d_k^*)$ is
defined by the equalities
\begin{equation}\label{12.02.2018.1}
d_1^\ast=\ldots =d_r^\ast=a +1, \hspace{1cm} d_{r+1}^\ast=
\ldots =d_k^\ast=a +2
\end{equation}
and $M=k a +(k-r)$, where} $0\leqslant r \leqslant
k-1$.\vspace{0.1cm}

{\bf Proof.} Explicitly, the proposition states that
$$
{M+l-e+m_e \choose M+l-e}\geqslant {M+l-e+m^*_e \choose M+l-e},
$$
where $m^*_e$ is calculated for the tuple $\underline{d}^*$. It is
easy to see that $m_e\geqslant m^*_e$, which proves the
proposition. Q.E.D. \vspace{0.1cm}

The second reduction simplifies the situation further, allowing us
to consider only the case when all degrees $d_i$ are
equal.\vspace{0.1cm}

{\bf Proposition 3.2.} {\it For the tuple
$\underline{d}^+=(d_1^+,\ldots , d_k^+)$ such that $d_1^+=\cdots
=d_k^+$, with $M^++k=|\underline{d}^+|$ and $M^+ \geqslant 8k\log
k - k$ the estimate
$$
\gamma(e,\underline{d}^+,l)\geqslant \frac{(M^+ -4k)(M^+
-5k)}{2}+M^+ +2k
$$
holds for all} $e=1,\dots, N_l^+=M^+-\mathop{\rm max}
\{[2\,\mathop{\rm log} k],l\}$.\vspace{0.1cm}

Let us show that Theorem 3.1 follows from Propositions 3.1 and
3.2.\vspace{0.1cm}

Indeed, by Proposition 3.1 it is sufficient to prove the
inequality
$$
\gamma(e,\underline{d}^*,l)\geqslant \frac{(M-5k)(M-6k)}{2}+M+k
$$
for $e=1,\dots, N_l$. Let us consider the tuple $\underline{d}^+$
with
$$
d^+_1=\cdots =d^+_k=a+1
$$
for the constant $a$ defined in Proposition 3.1. Set $M^+=ka$.
Obviously, $\gamma(e,\underline{d}^*,l)\geqslant
\gamma(e,\underline{d}^+,l)$ for $e=1,\dots,N^+_l$ as $M\geqslant
M^+$. If $N_l> N^+_l$, then for $i=0,\dots, N_l-N^+_l-1$ we have a
similar estimate ``from the other end'':
$$
\gamma(N_l-i,\underline{d}^*,l)\geqslant
\gamma(N^+_l-i,\underline{d}^+,l)
$$
(note that $N_l-N_l^+=M-M^+\leqslant k$). Therefore,
$$
\gamma(\underline{d}^*,l)=
\mathop{\rm min}\limits_{1\leqslant
e\leqslant N_l}\{\gamma(e,\underline{d}^*,l)\}\geqslant
\gamma(\underline{d}^+,l)=
\mathop{\rm min}\limits_{1\leqslant
e\leqslant N^+_l}\{\gamma(e,\underline{d}^+,l)\},
$$
and applying Proposition 3.2 and taking into account that
$M^+\geqslant M-k$, we get the claim of Theorem 3.1.\vspace{0.1cm}

The third reduction allows us to remove the integral parameter
$l\in \{0,1,\dots, k\}$. In order to simplify our notations, we
write $d_i$ for $d^+_i$, thus assuming that $d_1=\cdots =d_k=a+1$,
so that $M=ka$. We use the notation $\gamma(\underline{d},l)$ for
the minimum of the numbers $\gamma(e,\underline{d},l)$,
$e=1,\dots, N_l$, introduced above.\vspace{0.1cm}

{\bf Proposition 3.3.} {\it The following inequality holds:
$$
\gamma(\underline{d},l)\geqslant \gamma(\underline{d},0)
$$
for all} $l=0,1,\dots, k$.\vspace{0.1cm}

{\bf Proof.} Since for $l\geqslant 1$ we have $N_0>N_l$, it is
sufficient to compare the integers $\gamma(e,\underline{d},l)$ and
$\gamma(e,\underline{d},0)$ for the same values of $e=1,\dots,
N_l$. They are
$$
{M+l-e+m_e \choose M+l-e}\quad \mbox{and}
\quad {M-e+m_e \choose M-e},
$$
so the claim becomes obvious. Q.E.D.\vspace{0.1cm}

{\bf Remark 3.1.} We could as well do the third reduction as the
first one: show that the minimum of the integers
$\gamma(e,\underline{d},l)$ is attained for $l=0$ (which
corresponds to regular non-singular points of $V$), and after that
prove that the worst estimates correspond to the case
(\ref{12.02.2018.1}).\vspace{0.1cm}

The last (fourth) reduction makes the computations more compact.
Recall that now all degrees $d_i$ are equal to $a+1$. Introduce
the integer-valued function $\beta\colon \{2,\dots, a\}\to
{\mathbb Z}_+$ by the formula
$$
\beta(t)={k(a-t+1)+t \choose t}={kb(t)+t \choose t},
$$
where $b(t)=a-t+1$. Set also
$$
\alpha=\alpha(M,k)={a+1+[2 \log k]\choose a+1}.
$$

{\bf Proposition 3.4.} {\it The following estimate holds:}
$$
\gamma(\underline{d},0)\geqslant
\mathop{\rm min}\left\{\mathop{\rm min}\limits_{t\in\{2,\dots,a\}}
\{\beta(t)\},\alpha\right\}.
$$

{\bf Proof.} This follows immediately from the fact that for the
special tuple $\underline{d}$ of equal degrees
$$
m_{ki+1}=m_{ki+2}=\cdots =m_{ki+k}=i+2
$$
for $i=0,\dots, a-1$. Q.E.D.\vspace{0.1cm}

Therefore, the statement of Theorem 3.1 is implied by the
following facts. In the both propositions below we assume that $M
\geqslant 8k\log k - k$ and $k\geqslant 20$.\vspace{0.1cm}

{\bf Proposition 3.5.} {\it The minimum of the function $\beta(t)$
on the set $\{2,3,\dots, a\}$ is attained at $t=2$.}\vspace{0.1cm}

{\bf Proposition 3.6.} {\it The following inequality holds:}
$$
\alpha(M,k)\geqslant A(M,k)=\frac{(M-4k)(M-5k)}{2}+(M+2k).
$$

{\bf Remark 3.2.} The proof of Proposition 3.5 only requires $k
\geqslant 10$, it is Proposition 3.6 that requires $k \geqslant
20$, see a more detailed Remark 3.3. \vspace{0.1cm}

The rest of this section is a proof of the last two propositions,
which requires some (quite non-trivial) analytic
arguments.\vspace{0.3cm}

%%%%%%%%%%%%%%%%%%%%%%%%%%%%%%%%%%%%%%%%%%%%%%%%%%%%%%%%%%%%%%%%%%%%%%%%
%%%%%%%%%%%%%%%%%%%%%%%%%   subsection 3.3

{\bf 3.3. The Stirling formula.} The strategy of the proof of
Proposition 3.5 is as follows. Using the Stirling formula, we
construct a smooth function $\varepsilon\colon {\mathbb R}_+\to
{\mathbb R}$ such that $\varepsilon(t)\leqslant \beta(t)$ for
$t=2,\dots, a$ and $\varepsilon$ approximates $\beta$ with a good
precision. Then we show that the minimum of the function
$\varepsilon(t)$ on the interval $[2,a]$ occurs at one of the end
points, either for $t=2$ or for $t=a$. From this we deduce the
claim of Proposition 3.5.\vspace{0.1cm}

Recall that by the Stirling formula
$$
n!=\sqrt{2\pi n} n^n\exp(-n)\exp\left(\frac{\theta_n}{12n}\right)
$$
for some $\theta_n$ between $0$ and $1$. The integral parameter
$e$, enumerating the polynomials $g_e$, will never be used again
in this paper, so we use the symbol $e$ for the number $\exp(1)$.
Set
$$
\varepsilon(t)= \frac{\sqrt{2\pi}}{e^2}(kb(t)+t)^{(kb(t)+t+
\frac{1}{2})}(kb(t))^{-(kb(t)+ \frac{1}{2})}t^{-(t+\frac{1}{2})},
$$
by the Stirling formula $\beta(t)\geqslant \varepsilon
(t)$.\vspace{0.1cm}

{\bf Lemma 3.1.} {\it The smooth function $\varepsilon (t)$  for
$k\geqslant 3$ has only one critical point on the interval
$[2,a]$, which is a maximum, so that the minimum of that function
is attained at one of the end points.}\vspace{0.1cm}

{\bf Proof.} This is shown by demonstrating that\vspace{0.2cm}

(1) for $2 \leqslant t \leqslant \frac{M+k}{2k}$ the function
$\log \varepsilon (t)$ is strictly increasing,\vspace{0.2cm}

(2) for $\frac{M+1}{k+1} \leqslant t \leqslant a=\frac{M}{k}$ it
is strictly decreasing,\vspace{0.2cm}

(3) for $\frac{M+k}{2k} \leqslant t \leqslant \frac{M+1}{k+1}$ the
second derivative of $\log \varepsilon (t)$ is strictly negative
(this is where the maximum lies).\vspace{0.2cm}

The first derivative $\frac{d}{dt}{\log}\varepsilon(t)$ is equal
to
\begin{equation}\label{diff:1}
\frac{t^2-kb(t)^2}{2b(t)t(kb(t)+t)}-k \log \left(1+\frac{t}{kb(t)}
\right)+ \log \left( 1 + \frac{kb(t)}{t} \right),
\end{equation}
the second derivative $\frac{d^2}{dt^2}{\log}\varepsilon(t)$ is
given by the formula
\begin{equation}\label{diff:2}
\frac{1}{b(t)t}+\frac{(t^2-kb(t)^2)^2}{2b(t)^2t^2(kb(t)+t)^2}
+\frac{(k-1)(t^2-kb(t)^2)}{b(t)t(kb(t)+t)^2}-
\frac{k(t+b(t))^2}{tb(t)(kb(t)+t)}.
\end{equation}
We present the derivatives in these forms in order to use the
inequality
\begin{equation}
\label{ine:1}
\left|\frac{t^2-kb(t)^2}{2b(t)t(kb(t)+t)}
\right| \leqslant \frac{1}{2b(t)}.
\end{equation}

Now let us consider the domains (1)-(3) separately.\vspace{0.1cm}

(1) Assume that $2 \leqslant t \leqslant \frac{M+k}{2k}$. Note
that on this interval $b(t)\geqslant 2$ so that
$$
\left|\frac{t^2-kb(t)^2}{2b(t)t(kb(t)+t)} \right|
\leqslant \frac{1}{4}.
$$
The last term in the expression (\ref{diff:1}) can be estimated as
$$
\log \left( 1 + \frac{kb(t)}{t} \right) \geqslant \log(1+k)
\geqslant \log 21 >3,
$$
since on the interval $[2, \frac{M+k}{2k}]$ we have $t \leqslant
b(t)$. Finally, for the second term in (\ref{diff:1}) we get
$$
-k \log \left(1+\frac{t}{kb(t)} \right) \geqslant -\frac{t}{b(t)}
\geqslant -1.
$$
Combining these estimates, we obtain the inequality
$$
\left.\frac{d}{dt} \log \varepsilon(t)\right|_{2 \leqslant
t\leqslant \frac{M+k}{2k}} \geqslant -\frac{5}{4} + 3 > 0,
$$
so that indeed $\varepsilon (t)$ is increasing on the interval
under consideration.\vspace{0.1cm}

(2) Assume now that $\frac{M+1}{k+1} \leqslant t \leqslant
\frac{M}{k}$. Here $t \geqslant kb(t)\geqslant k$. First of all,
we have the inequality
$$
\left|\frac{t^2-kb(t)^2}{2b(t)t(kb(t)+t)} \right| \leqslant \frac{1}{2}.
$$
For the other two terms in the expression (\ref{diff:1}) we get
the estimates
$$
-k \log \left(1+\frac{t}{kb(t)} \right) \leqslant -k \log 2
$$
and
$$
\log \left( 1 + \frac{kb(t)}{t} \right) \leqslant \log 2.
$$
Combining these inequalities, we see that
$$
\frac{d}{dt} \log \varepsilon(t)\leqslant \frac{1}{2} - (k-1) \log
2 < 0
$$
for $t\in [\frac{M+1}{k+1},\frac{M}{k}]$ as we claimed
above.\vspace{0.1cm}

(3) Finally, assume that $\frac{M+k}{2k} \leqslant t \leqslant
\frac{M+1}{k+1}$. On this interval $b(t) \leqslant t \leqslant
kb(t)$. Let us show that the second derivative (\ref{diff:2}) is
negative. Using again the inequality (\ref{ine:1}), we get that
$\frac{d^2}{dt^2} \log \varepsilon(t)$ on the interval under
consideration is not higher than
$$
\frac{1}{b(t)t}+\frac{1}{2b(t)^2}+
\frac{(k-1)}{b(t)(kb(t)+t)}-\frac{k(t+b(t))^2}{tb(t)(kb(t)+t)}=
$$
$$
=\frac{t^2+kb(t)(-2t^2-4b(t)t+3t-2b(t)^2+2b(t))}{2tb(t)^2(kb(t)+t)}.
$$
Elementary computations, together with the inequality $t\leqslant
a$, show that the expression in brackets in the numerator is not
higher than $-2a^2-a$. Therefore, the whole numerator is not
higher than
$$
t^2-kb(t)(2a^2+a)\leqslant kb(t)(t-2a^2-a)<0.
$$
We have shown that $\frac{d^2}{dt^2} \log \varepsilon(t)<0$ for
$t\in [\frac{M+k}{2k}, \frac{M+1}{k+1}]$. This completes the proof
of Lemma 3.1. Q.E.D.\vspace{0.3cm}

%%%%%%%%%%%%%%%%%%%%%%%%%%%%%%%%%%%%%%%%%%%%%%%%%%%%%%%%%%%%%%%%%%%%%%%%%%%
%%%%%%%%%%%%%%%%%%%%%%%%%%%%%%%%%%%%%%%%%%%%%%   subsection 3.4

{\bf 3.4. Proof of Proposition 3.5.} In view of the inequality
$\varepsilon(t)\leqslant\beta(t)$ and Lemma 3.1, Proposition 3.5
follows from the two lemmas stated below.\vspace{0.1cm}

{\bf Lemma 3.2.} {\it The inequality
$\beta(2)\leqslant\varepsilon(3)$ holds.}\vspace{0.1cm}

{\bf Lemma 3.3.} {\it The inequality
$\beta(2)\leqslant\varepsilon(a)$ holds.}\vspace{0.1cm}

{\bf Proof of Lemma 3.2.} We need to estimate the error in
Stirling's approximation, in order to be able to use $\beta(3)$
instead of $\varepsilon(3)$. The number $\beta(3)$ is a polynomial
in $M,k$, which makes the task easier. From the Stirling formula
we get:
$$
1.126\cdot \varepsilon(3) \leqslant \beta(3)
\leqslant 1.132\cdot \varepsilon (3).
$$
The inequality of the lemma will follow if it is shown that $1.14
\cdot \beta(2) < \beta(3)$ and this is equivalent to the
inequality $G_1(M,k)=6(\beta(3)-1.14 \cdot \beta(2))> 0$. Here
$G_1(M,k)$ is given explicitly by the expression
$$
M^3+M^2(2.58-6k)+M(12k^2-17.16k+0.74)-8k^3+20.58k^2-11.74k-0.84.
$$
It is easy to check that $G_1(8k\log k,k)$ and the partial
derivative $\frac{\partial}{\partial M}G_1(M,k)$ are both positive
for $k \geqslant 20$ and $M\geqslant 8k\log k$. This completes the
proof. Q.E.D.\vspace{0.1cm}

{\bf Proof of Lemma 3.3.} The claim of the lemma is equivalent to
the inequality
$$
G_2(M,k)= \log \varepsilon(a)- \log \beta(2)\geqslant 0.
$$
A direct calculation gives $G_2(160 \log 20 - 20, 20)>0$. Set
$$
G_3(t)=\frac{d}{d t} G_2(8t\log t - t, t).
$$

{\bf Lemma 3.4.} $G_3(t)>0$ {\it for} $t\geqslant
20$.\vspace{0.1cm}

{\bf Proof.} Explicitly,
$$
G_3(t) = \log \left(1+\frac{8\log t -1}{t}\right)+ \frac{8}{t}
\log \left(1+\frac{t}{8\log t - 1}\right)-
\frac{1}{2t}+H_1(t)+H_2(t),
$$
where
$$
H_1(t)=\left(\frac{8}{t}+1 \right) \frac{8\log t+t - 0.5}{8\log
t+t - 1}-\left(\frac{8}{t} \right)\frac{8\log t - 0.5}{8\log t - 1
}-1,
$$
$$
H_2(t)=-(8\log t+6)\left(\frac{1}{8t\log t-2t+2}+ \frac{1}{8t\log
t-2t+1} \right).
$$
Using the power series expansion of $\log (1+x)$, we obtain the
inequality
$$
G_3(t)>-\frac{1}{2t}+\frac{8\log t - 1}{t}- \frac{(8 \log t -
1)^2}{2t^2}+\frac{8}{t} \log \left(1+\frac{t}{8\log t - 1}\right)
+H_1(t)+H_2(t).
$$
For $t \geqslant 20$ then $H_1(t)\geqslant 0$ and $H_2(t)
\geqslant -\frac{4}{t}$, which can be checked directly. This gives
the inequality
$$
G_3(t)>\frac{16\log t- 11}{2t}-\frac{(8 \log t - 1)^2}{2t^2}+
\frac{8}{t} \log \left(1+\frac{t}{8\log t - 1}\right).
$$
The right hand side of the last inequality is higher than
$$
\frac{1}{2t^2}(16t\log t - 11 t -64\,(\log t)^2 + \log t - 1),
$$
which is positive for $t\geqslant 20$. Q.E.D. for Lemma
3.4.\vspace{0.1cm}

We conclude that $G_2(8t\log t - t, t)>0$ for $t\geqslant 20$. The
claim of Lemma 3.3 will be proven if we show that for $k \geqslant
20$ and $M\geqslant 8k\log k - k$ the function $G_2(M,k)$  is an
increasing function of $M$. Set
$$
G_4(s,t)=\frac{\partial}{\partial s} G_2(s,t).
$$

{\bf Lemma 3.5.} $G_4(s,t)>0$ {\it for} $t\geqslant 20, s\geqslant
8t\log t - t$.\vspace{0.1cm}

{\bf Proof.} Explicitly,
$$
G_4(s,t)=\frac{1}{t} \log \left(1+\frac{t^2}{s} \right)-
\frac{t^2}{2s(t^2+s)}-\frac{2s+3-2t}{s^2+(3-2t)s+t^2-3t+2}.
$$
First we consider the case when $s\leqslant t^2$ and get
$$
\left.G_4\right|_{s \leqslant t^2} \geqslant \frac{1}{t} \log
2-\frac{t^2}{2s(t^2+s)}- \frac{2s+3-2t}{s^2+(3-2t)s+t^2-3t+2}.
$$
It is easy to see that the minimum of the right hand side occurs
when $s=8t\log t - t$ is the smallest possible, so that for
$s\leqslant t^2$ the function $G_4(s,t)$ is bounded from below by
the expression
$$
\frac{1}{t} \log 2 - \frac{1}{(16 \log t - 2)(t+8 \log t - 1)}-
$$
$$
- \frac{16t\log t+3-4t}{t^2(8 \log t - 1)^2+(3-2t)(8t\log t -
t)+t^2-3t+2},
$$
which is positive for $t \geqslant 20$.\vspace{0.1cm}

Now let us consider the region $s \geqslant t^2$. Here we get
$$
G_4(s,t) \geqslant \frac{t}{s}-\frac{t^3}{2s^2}-\frac{t^2}{2s(t^2+s)}-
\frac{2s+3-2t}{s^2+(3-2t)s+t^2-3t+2}.
$$
A direct check shows that for $t \geqslant 20$ the expression in
the right hand side is positive. Q.E.D. for Lemmas 3.5, 3.3 and
Proposition 3.5. \vspace{0.3cm}

%%%%%%%%%%%%%%%%%%%%%%%%%%%%%%%%%%%%%%%%%%%%%%%%%%%%%%%%%%%%%%%%%%%%%%%%
%%%%%%%%%%%%%%%%%%%%%%%%%%%%%   subsection 3.5

{\bf 3.5. Proof of Proposition 3.6.} This proof is obtained in the
same way as that of Proposition 3.5 and we only point out the main
steps of the computations, leaving the details to the reader. In
order to prove the inequality $\alpha(M,k)\geqslant A(M,k)$, we
use the Stirling approximation of $\alpha(M,k)$. Namely, we
introduce the function $G_5(s,t,r)$ of three real variables by the
formula
$$
G_5(s,t,r)= \left(\frac{s}{t}+r+\frac{3}{2}\right) \log
\left(\frac{s}{t}+r+1\right)- \left(r +\frac{1}{2}\right) \log r-
$$
$$
-\left(\frac{s}{t}+\frac{3}{2}\right) \log
\left(\frac{s}{t}+1\right)+ \log
\left(\frac{\sqrt{2\pi}}{e^2}\right)- \log A(s,t).
$$
By the Stirling approximation, Proposition 3.6 follows from the inequality
$$
G_5(M, k, [2 \log k])\geqslant 0.
$$
It is easy to see that
$$
G_5(M, k, [2 \log k])\geqslant G_5(M, k, 2\log k - 1),
$$
so we set $G_6(s,t)=G_5(s, t, 2\log t-1)$ and prove the inequality
$$
G_6(s,t)\geqslant 0
$$
for $s\geqslant 8t\log t - t$, $t\geqslant 20$. First of all,
explicit computations show that
$$
G_6(8t\log t - t,t)\geqslant 0
$$
for $t\geqslant 20$. Set
$$
G_7(s,t)=\frac{\partial}{\partial s} G_6(s,t).
$$
It remains to show that for $s\geqslant 8t\log t - t$ and
$t\geqslant 20$ the inequality $G_7(s,t)\geqslant 0$ holds.
Explicitly, $G_7(s,t)$ is given by the expression
$$
\frac{1}{t} \log \left(1+\frac{2\log t-1}{\frac{s}{t}+1} \right)
-\frac{2\log t-1}{2t(\frac{s}{t}+1)(\frac{s}{t}+2\log t)}
-\frac{2s-9t+2}{s^2-9ts+2s+20t^2+4t}.
$$
Now the inequality $G_7(s,t)\geqslant 0$ is obtained by tedious
but straightforward computations, using the estimate
$$
\frac{1}{t}{\log}\left(1+\frac{2\log t-1}{\frac{s}{t}+1} \right)
>\frac{2\log t-1}{t(\frac{s}{t}+1)}- \frac{(2\log
t-1)^2}{2t(\frac{s}{t}+1)^2}.
$$
The details are left to the reader. Q.E.D. for Proposition 3.6 and
Theorem 3.1.\vspace{0.1cm}

{\bf Remark 3.3.} (i) It is clear from the computations presented
in this section and in the proof of Lemma 1.3 that in certain
parts of our arguments we need much weaker lower bounds for $k$.
For instance, Lemma 3.2 requires only that $k\geqslant 5$ and
Lemmas 3.3 and 3.4 require only that $k\geqslant 10$. In the last
part of Subsection 3.5, for the inequality $G_7(s,t)\geqslant 0$
only the assumption $t\geqslant 10$ is needed. However, the
inequality
$$
G_6(8t\log t - t,t)\geqslant 0
$$
requires that $t\geqslant 20$. In order for the whole argument to
work, we have to select the strongest restriction.\vspace{0.1cm}

(ii) One more paper \cite{Zhuang2018} on birational superrigidity
of non-singular Fano complete intersections of index one was put
on the archive when our paper was finalized. We add it on the
reference list.

%%%%%%%%%%%%%%%%%%%%%%%%%%%%%%%%%%%%%%%%%%%%%%%%%%%%%%%%%%%%%%%%%%%%%%%%%%%%%%%%
%%%%%%%%%%%%%%%%%%%%%%%%%%%%%%%%%%%%%%%%%%%%%%%%%%%%%%%%%%%%%%%%%%%%%%%%%%%%%%%%
%%%%%%%%%%%%%%%%%%%%%%%%%%%%%%%%%%%%%%%%%%%%%%%%%%%%%%%%%%%%%%%%%%%%%%%%%%%%%%%%

\begin{flushleft}
Department of Mathematical Sciences,\\
The University of Liverpool
\end{flushleft}

\noindent{\it pukh@liverpool.ac.uk, sgdevans@liverpool.ac.uk}

\end{document}